\documentclass{amsart}
\usepackage{amsmath}
\usepackage{amssymb}
\usepackage{amsthm}
\usepackage{amsfonts}

\def\H{{\mathcal H}}
\def\C{\mathbb{C}}
\def\R{\mathbb{R}}

\def\End{\textup{End}}
\def\dim{\textup{dim }}
\def\X{{\mathcal{X}}}

\newcommand{\G}[1]{\textup{GL}_{#1}(\C)}

\newcommand{\vect}[2]{{#1}_{{#2} {\bullet}}}
\newcommand{\colvect}[2]{{#1}_{{\bullet}{#2}}}

\theoremstyle{plain}
\newtheorem{theorem}{Theorem}
\newtheorem{lemma}[theorem]{Lemma}
\newtheorem{proposition}[theorem]{Proposition}
\newtheorem{corollary}[theorem]{Corollary}
\newtheorem{question}{Question}
\numberwithin{theorem}{section}

\theoremstyle{definition}
\newtheorem{definition}[theorem]{Definition}

\numberwithin{equation}{section}

\title{Linear conditions imposed on flag varieties}
\author{Julianna S. Tymoczko} 

\address{Department of Mathematics, University of Michigan, 2074 East Hall, Ann Arbor, MI 48109-1109}
\email{tymoczko@umich.edu}

\begin{document}
\begin{abstract}
We study subvarieties of the flag variety called Hessenberg varieties,
 defined by certain linear conditions.  
These subvarieties arise naturally in 
applications including geometric representation theory, number theory, and 
numerical analysis.  We describe completely the homology of Hessenberg
varieties over $\G{n}$ and show that they have no odd-dimensional
homology.  We provide an explicit geometric construction which
partitions each Hessenberg variety into pieces homeomorphic to 
affine space.  We  characterize these
affine pieces by fillings of Young tableaux 
and show that the dimension of the affine piece can be computed 
by combinatorial rules generalizing the Eulerian numbers.  We give an 
equivalent formulation of this result in terms
of roots.  We conclude with a section on open questions.
\end{abstract}

\maketitle

\section{Introduction}

The 
full flag variety over $\G{n}$ is the collection of nested complex vector
spaces $V_1 \subsetneq V_2 \subsetneq \cdots \subsetneq V_n = \C^n$
where 
$V_i$ is $i$-dimensional.  Given a linear operator $X$ on $\C^n$, the
set of flags that are stabilized by $X$---that is, flags 
$V_1 \subsetneq \cdots \subsetneq V_n$ such that $X V_i \subseteq V_i$
for each $i$---is an important subvariety of the full flag variety called
the Springer-Grothendieck fiber.  Geometric representation theorists
 use this subvariety to construct the irreducible representations
of the symmetric group
(\cite[section 3.6]{CG} has background and references).

More generally, fix any nondecreasing 
function $h: \{1, 2, \ldots, n\} \rightarrow \{1, 2, \ldots,
n\}$ such that $h(i) \geq i$ for each $i$, and consider the flags
\[\H(X,h)=\{ \textup{flags } V_1 \subseteq \cdots \subseteq V_n
  \textup{ such that }
   X V_i \subseteq V_{h(i)} \textup{ for each }i\}.\]
The subvariety $\H(X,h)$ is called a Hessenberg variety, and the map $h$ 
is a Hessenberg function.

For example, consider the set of flags with $X V_i \subseteq V_{i+1}$
whenever $i$ is less than $n$.  
This parametrizes the bases that put the operator $X$
into Hessenberg form, a form used in 
a common algorithm to compute eigenvalues (see \cite{dMS} for more
about the QR algorithm).  The natural generalization presented here 
was defined in \cite{dMPS}.

Our main theorem explicitly partitions
each Hessenberg variety into affine spaces satisfying 
weak closure rules.  This decomposition is a paving and 
is the intersection of $\H(X,h)$ with 
a special Bruhat decomposition of the flag variety.
Pavings give the homology of $\H(X,h)$, and hence  a
combinatorial description of its Betti numbers.  We  conclude that
Hessenberg varieties have no odd-dimensional homology.  

For notational convenience, we give the main result here in the case
when $X=N$ is nilpotent.  Theorems \ref{main theorem: bundle form} and 
\ref{theorem: tableau form}
have the result for general $X$ 
in two different forms.  If the nilpotent operator $N$  has
Jordan blocks  of size $d_1 \geq \ldots \geq d_k$, then associate
to it the Young diagram $\lambda_N$ with row lengths $d_1 \geq 
\ldots \geq d_k$.  Our Young diagrams are left-aligned and 
top-aligned.  For example, Figure \ref{defn tableaux} shows a nilpotent
with Jordan blocks of dimension $3$ and $1$ and the associated 
Young diagram.

\begin{figure}[h]
\[\mbox{\small $\left( \begin{array}{cccc}
0 & 1 & 0 & 0 \\
0 & 0 & 1 & 0 \\
0 & 0 & 0 & 0 \\
0 & 0 & 0 & 0
\end{array} \right)$}
\longleftrightarrow 
\textup{\begin{tabular}{|c|c|c|}
\cline{1-3} & & \\
\cline{1-3} & \multicolumn{2}{c}{} \\
\cline{1-1} \multicolumn{3}{c}{} \vspace{-1.3em} \end{tabular}} \]
\caption{The Young diagram corresponding to one nilpotent operator}\label{defn tableaux}
\end{figure}

The cells of the paving are indexed by Young tableaux 
that are filled with the numbers 
from $1$ to $n$ without repetition.  Each tableau defines a permutation  
$w$ of $n$ letters for which 
$w^{-1}(k)$ is the number of boxes to the left of or below
the box filled by $k$ (including the box itself).

\begin{figure}[h]
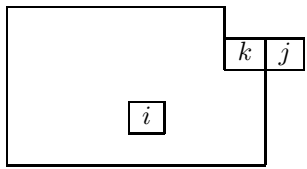

\begin{tabular}{|c|c|c|c|c|}
\cline{1-3} \multicolumn{3}{|c|}{$\hspace{1in} $} & 
    \multicolumn{2}{c}{} \\ 
\cline{4-5} \multicolumn{3}{|c|}{} & $k$ & $j$ \\
\cline{4-5} \multicolumn{3}{|c}{} & \multicolumn{1}{c|}{} &
    \multicolumn{1}{c}{} \\
\cline{2-2} \multicolumn{1}{|c|}{$\hspace{.5in}$} & \multicolumn{1}{|c|}{$i$}
   & \multicolumn{2}{|c|}{} & \multicolumn{1}{c}{} \\
\cline{2-2} \multicolumn{3}{|c}{} & \multicolumn{1}{c|}{} &
    \multicolumn{1}{c}{} \\
\cline{1-4} \multicolumn{5}{c}{}
\end{tabular} \vspace{-.15in}
\caption{Configuration of triples} \label{picture of conditions}
\end{figure}

\begin{theorem} \label{nilpotent: tableau form}
Fix a nilpotent $N$.  
The Hessenberg variety $\H(N,h)$ is paved by affines.  Each nonempty
cell corresponds to a unique  
filling of $\lambda_N$ in which 
\begin{tabular}{|c|c|} \cline{1-2} $k$ & $j$ \\ \hline
\end{tabular} occurs only if $k \leq h(j)$.  This correspondance is
a bijection.  The
dimension of a nonempty cell is the number of pairs $i$,$k$
such that 
\begin{enumerate}
\item $i$ is below or anywhere to the left of $k$ (see Figure 
\ref{picture of conditions}), 
\item $k < i$, and
\item \label{inequality} 
if there is a box immediately to the right of $k$ that is filled by
$j$ then $i \leq h(j)$.
\end{enumerate}
\end{theorem}


This result extends N.~Spaltenstein's description of the Springer
fibers' components, the case when $h(i)=i$ \cite{Sp}.
In particular, it can be used to give a new proof that the
rank of each irreducible representation of the symmetric group
is the number of standard fillings of its Young tableau.
It also partially extends 
the work of F.~de Mari, C.~Procesi, and M.~Shayman paving 
Hessenberg varieties by affines
when $X$ is regular semisimple \cite{dMPS}, and of
C.~de Concini, G.~Lusztig, and C.~Procesi paving Springer fibers 
by affines \cite{dCLP}.  
Our methods are different from theirs
though similar in spirit to Spaltenstein's or to those in \cite{KnM}.  
B.~Kostant used a different Bruhat decomposition to pave
one Hessenberg variety when $X$ is regular nilpotent, 
giving a geometric construction of the quantum cohomology of the
flag variety \cite{K}.  According to personal communications
\cite{C} and announcements \cite[Theorem 3]{BC}, D.~Peterson has 
other uncirculated results studying Hessenberg varieties when $X$ is
regular nilpotent. Our methods 
do not use torus actions, as there is no obvious torus action
for general $X$.  Rather than using one-dimensional deformations as
in \cite{V} or restricting to intersections with 
codimension-one Schubert varieties as
in \cite{So}, our approach makes fewer deformations of higher dimension
in each Schubert cell.

Our proof begins by describing $\H(X,h)$ in terms of matrices $g$
for which $g^{-1}Xg$ is zero in fixed coordinates, and then 
reducing to the case when $g=u$ is upper-triangular.  
The entries of 
the matrices $u^{-1}Xu$ need not be linear nor affine 
functions of the entries of $u$.  However, the entries of
the $i^{th}$ row of $u^{-1}Xu$ are affine functions of the $i^{th}$ row of
$u$.  For instance, when $X$ is nilpotent with a single Jordan block
its conjugate by an upper-triangular $u$ is
\[ u^{-1} \left( \begin{array}{cccc} 0 & 1 & 0 & 0 \\ 0 & 0 & 1 &0 
\\ 0 & 0 & 0 & 1 \\ 0 & 0 & 0 & 0
\end{array} \right) u =
\left( \begin{array}{cccc} 0 & 1 & u_{23}-u_{12} & u_{24} - 
    u_{12}(u_{34}-u_{23}) - u_{13} \\ 0 & 0 & 1 & u_{34}-u_{23}
\\ 0 & 0 & 0 & 1 \\ 0 & 0 & 0 & 0
\end{array} \right).\]
Section \ref{paving definitions} 
has the necessary background on the Bruhat decomposition and pavings.
In Section \ref{row section} 
we partition the upper-triangular matrices into subgroups
called {\em rows}
and show that conjugation by a row is 
an affine transformation of the row's entries.  

In our example, the functions of $u_{1j}$ in the first row have the same rank
regardless of the other $u_{ij}$.  This is  true if $X$ is
in {\em highest form}, defined for any linear operator in  
Section \ref{highest form section}.  
Section \ref{paving by affines section} 
has the paper's key lemma.  That lemma is one step in the main theorem
of Section \ref{main theorems}, which proves 
that each cell of a Bruhat decomposition 
intersects $\H(X,H)$ in an iterated tower of affine fiber bundles.
The main theorem is described using tableaux
 in Section \ref{main theorem: tableau form} and 
using roots in Section \ref{main theorem: root form}.
Section \ref{open questions} has open
questions and conjectures about Hessenberg varieties, including whether
they are pure dimensional and how many components they have.

This work was partially supported by an NDSE graduate fellowship and
was part of the author's doctoral dissertation.  I thank Emina 
Alibegovi\'{c},
Jared Anderson, Henry Cohn, William Fulton, Gil Kalai, David Kazhdan,
Robert Lazarsfeld, David Nadler, Arun Ram, Eric Sommers, and the  
anonymous referee for valuable comments.  I am especially 
grateful to my advisor, Robert MacPherson.

\section{Pavings and the Bruhat decomposition} \label{paving definitions}

In this section we describe a
classical partition of the flag variety called the Bruhat decomposition.
We also precisely define pavings, the special 
partitions of a variety used in this paper, sometimes called cellular
decompositions.  

\begin{definition} A paving of the variety $\X$ is an ordered partition
$\X = \coprod_{i=0}^{\infty} \X_i$ so that each finite union $\coprod_{i=0}^j
\X_i$ is Zariski-closed in $\X$.  
If in addition each $\X_i$ is homeomorphic to
affine space $\R^{d_i}$ then $\coprod_{i=0}^{\infty} \X_i$ 
is a paving by affines.
\end{definition}

\begin{figure}[h]
\begin{picture}(250,40)(0,-20)
\put(65,0){\circle*{5}}
\put(85,0){\circle{40}}
\put(125,0){\circle{40}}
\put(165,0){\circle{40}}
\put(85,0){\oval(40,20)[b]}
\put(125,0){\oval(40,20)[b]}
\put(165,0){\oval(40,20)[b]}
\end{picture}
\caption{A Space Paved by Four Cells} \label{paving example}
\end{figure}
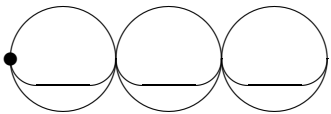

Our pavings have a finite number of pieces.  We
call the $\X_i$ cells.  Figure \ref{paving example}
shows three spheres glued successively at a point like a string of
beads.  It is paved by four affine cells: 
the marked point and each $S^2$ without its leftmost point.  
The closure of a cell need not cover the cells it intersects, 
as it must in a CW-decomposition.

Pavings by affines determine Betti numbers \cite[19.1.11]{F}:

\begin{proposition} \label{betti of pavings}
Let $\X=\coprod \X_i$ be a paving by a finite number of affines $\X_i$
with each $\X_i$ homeomorphic to $\C^{d_i}$.  Then the 
nonzero cohomology groups of $\X$ are 
\[H^{k}(\X) = \bigoplus_{i \textup{ such that }  2d_i=k} \mathbb{Z}.\]
\end{proposition}

The full flag variety has a well-known paving by
affines called the Bruhat decomposition.  Recall that 
the flag $V_1 \subseteq \cdots \subseteq V_n$ is determined 
by any matrix $g$ 
whose first $i$ column vectors generate the $i^{th}$ vector space $V_i$.  
The flag corresponding to $g$ is denoted $[g]$.

The next definition parametrizes the cells of this paving 
\cite[section 28.4]{H}.  

\begin{definition}
Let $w$ be a permutation matrix.  The group $U_w$ of upper-triangular matrices 
associated to $w$ is defined as $U_w = \{u: u \in U, w^{-1}uw \textup{ is
lower-triangular}\}$.
\end{definition}

We now state a classical result 
in the language of this paper.   Write $e_i$ for the
basis vector of $\C^n$ which has one in the $i^{th}$ position and zero
otherwise.  The permutation matrix $w$ corresponds to the
permutation of $\{1,2,\ldots,n\}$ given by $e_i w = e_{w(i)}$.

\begin{proposition} \label{Bruhat decomposition}
The flag variety is paved by affines $\coprod_{w \in S_n} C_w$.  
The Schubert cell $C_w$ is the set of flags $[U_w w]$, which is 
homeomorphic to $U_w w$ and has 
dimension $|\{(i,j): 1 \leq i<j \leq n, w(i) > w(j)\}|$.
\end{proposition}

\begin{proof}
The Schubert cells are described in \cite[section 28.3]{H}.
The $U_w$ parametrize the cells by \cite[section 28.4]{H}. 
The cells form a paving by \cite[section 2.10]{BL}.
\end{proof}

The matrix description of the flag variety gives a different formulation of
the definition of Hessenberg varieties.  

\begin{definition} \label{Hessenberg space}
The Hessenberg space $H$ associated to $h$ is the linear subspace of 
matrices $X$ whose $(i,j)^{th}$ entry $X_{ij}=0$ if $i>h(j)$.
\end{definition}

Section \ref{main theorem: root form} has an intrinsic
definition of Hessenberg spaces from \cite{dMPS}.
The next proposition relates the linear subspace $H$ to the
function $h$.  Its proof is immediate from $w^{-1}E_{jk}w = E_{w(j),w(k)}$,
where $E_{jk}$ is the matrix basis unit with $1$ in its $(j,k)$ entry and
zero everywhere else.

\begin{proposition} \label{algebraic to combinatorial}
The matrix basis unit $E_{jk} \in wHw^{-1}$ if and only if $w(j) 
\leq h(w(k))$. 
\end{proposition}

An alternate definition of Hessenberg varieties first given in \cite{dMPS}
is
\[\H(X,H) = \{\textup{flags }[g]: g^{-1}Xg \in H\} = \H(X,h).\]

Conjugation by 
$g \in GL_n(\C)$ is a homeomorphism of Hessenberg varieties in two 
ways.

\begin{proposition}
Fix $X$ and $H$ and $g_0 \in GL_n(\C)$.  
The Hessenberg variety $\H(g_0^{-1}Xg_0,H)$ is
homeomorphic to $\H(X,H)$.
\end{proposition}

\begin{proof}  Using associativity gives
$\H(g_0^{-1}Xg_0,H) = g_0^{-1} \H(X,H)$.
Multiplication is an automorphism of flags so this
is homeomorphic to $\H(X,H)$.  
\end{proof}

\begin{proposition}
Fix a matrix $X$, a Hessenberg space $H$, and $g_0 \in GL_n(\C)$.  
The Hessenberg variety $\H(g_0^{-1}Xg_0,g_0^{-1}Hg_0)$ is homeomorphic
to $\H(X,H)$.
\end{proposition}

\begin{proof} By definition,
$\H(g_0^{-1}Xg_0,g_0^{-1}Hg_0) 
= \{\textup{flags }[g_0^{-1} g g_0]: g^{-1}Xg \in H\}$.
Conjugation is an automorphism of flags so
this is homeomorphic to $\H(X,H)$.
\end{proof}

These show that the topology and geometry
of an arbitrary Hessenberg variety $\H(X,H)$ are the same as when
$X$, $H$, and the underlying basis are in fixed relative position.
In what follows, we assume that $X$ and $H$ are in fixed conjugacy classes
without further comment.

\section{Rows of upper-triangular matrices} \label{row section}

This section describes a decomposition of the upper-triangular
invertible matrices into subgroups called rows and shows how rows
act on arbitrary matrices.  A similar partition is used implicitly
in \cite[section 2.C]{Ste} and
in \cite[section 3]{CP}.

Unless otherwise stated all
matrices are $n \times n$ with complex coefficients.
We use $X$ to denote an arbitrary matrix, $N$ to denote 
a nilpotent upper-triangular
matrix, and $S$ to denote a diagonal matrix.  Write $U$
for the group of upper-triangular matrices with ones on the diagonal.
Let $X_{jk}$ be the $(j,k)^{th}$ entry of the matrix $X$.  

\begin{definition}
The $i^{th}$ row $U_i$ is the subgroup 
$U_i = \left\{u \in U: u_{jk}=0 \textup{ if }
  j \neq i,k\right\}$.
\end{definition}

We distinguish the rows $U_i$ from the Schubert cell subgroups $U_w$ 
by  subscripts: $i$, $j$, $k$ always denote an integer, while $w$ 
always denotes a permutation matrix.  Note that
$U_i \cap U_j$ is the identity if $i \neq j$.  The 
rows generate all of $U$ because each row is a product of one-parameter
subgroups, as in \cite[Proposition 28.1]{H}.

\begin{proposition} \label{rows factor U}
The group $U$ factors uniquely as $U=U_{n-1} U_{n-2} \cdots U_1$.
\end{proposition}

This result together with Proposition \ref{Bruhat decomposition} 
shows that representatives for each Schubert cell factor uniquely as 
$(U_w \cap U_{n-1}) (U_w \cap U_{n-2}) \cdots (U_w \cap U_1)w$.

We use rows because of their group structure, given
next.  Its proof is immediate.

\begin{proposition} \label{row structure}
$U_i$ is naturally isomorphic to the additive group $\C^{n-i}$.
If $u$ and $v$ are elements of $U_i$
then $(uv)_{ik}=u_{ik}+v_{ik}$ for each $k>i$.  
In particular, the entries of the inverse $u^{-1}$ are given by 
$\left(u^{-1}\right)_{ik}=-u_{ik}$ for each $k > i$.
\end{proposition}

The group $U_i$ acts on a matrix $X$ by left-multiplication,
right-multiplication, or conjugation.  In each
case most of the rows of $X$ are preserved, as the following makes precise.

\begin{proposition} \label{row conjugation}
Fix $u$ in $U_i$.
\begin{enumerate} 
\item \label{row difference} $(uX)_{jk} = X_{jk}$ 
  except possibly when $j=i$.
\item \label{column difference} $(Xu)_{jk}=X_{jk}$ 
  except possibly in rows $j$ for which $X_{ji}$ is nonzero.
\item \label{conjugate difference} 
  If $X$ is upper-triangular then $\left(u^{-1}Xu \right)_{jk} = X_{jk}$ 
  except possibly when $j \leq i$.
\end{enumerate}
\end{proposition}

\begin{proof}
The first two parts restate matrix multiplication.

Since $X$ is upper triangular
the product $(Xu)_{jk}=X_{jk}$ except perhaps in a row $j$ with $j \leq i$
by Part \ref{column difference}.
By Part \ref{row difference} the product $\left(u^{-1}Xu\right)_{jk} = 
(Xu)_{jk}$ except perhaps when $j=i$.  Thus $\left(u^{-1}Xu \right)_{jk}=
X_{jk}$ whenever $j > i$.
\end{proof}

Denote the $i^{th}$ row vector of $X$ by $\vect{X}{i}$.  Let $X=S+N$ be
upper-triangular.
The next result shows that the $i^{th}$ row of $u^{-1}Xu$
is the image under an affine transformation of the $i^{th}$ row of $u$, namely
the translation of a linear map on $\vect{u}{i}$.

\begin{proposition}\label{affine transformation}
The map $\vect{u}{i} \mapsto \vect{u^{-1}(S+N)u}{i}$ is an affine
transformation of the entries of $\vect{u}{i}$.  Explicitly,
\[ \vect{\left( u^{-1}(S+N)u \right)}{i} = S_{ii} 
   \vect{u}{i} + \vect{\left(u^{-1}\right)}{i} (S+N).\]
\end{proposition}

\begin{proof}
We prove this by comparing the $k^{th}$ entry of each vector.  Note that
\[\begin{array}{rl} 
\displaystyle{\left( u^{-1}(S+N)u \right)_{ik}} & \displaystyle{= 
    \sum_{j=1}^n \left( u^{-1} \right)_{ij}
     \left( (S+N)u \right)_{jk}} \\
     &\displaystyle{= \sum_{j=1}^n \left( u^{-1} \right)_{ij} 
       (S_{ji}+N_{ji}) u_{ik}
        + \sum_{j=1}^n \left( u^{-1} \right)_{ij}
         (S_{jk}+N_{jk}) u_{kk}.}
\end{array}\]
The first sum simplifies to $(u^{-1})_{ii}(S_{ii} + N_{ii})u_{ik}$
because if $i>j$ then $(u^{-1})_{ij}=0$ and if $i<j$ both $S_{ji}$
and $N_{ji}$ vanish.  Since $N_{ii}=0$ and $(u^{-1})_{ii}=1$ this
is $S_{ii}u_{ik}$.

The second sum is  the $k^{th}$ entry of $\vect{(u^{-1})}{i} (S+N)$
by definition.
\end{proof}

\section{Highest forms of linear operators} \label{highest form section}

This section introduces one of the
main tools of our proof: the {\em highest form} 
for linear operators.  We first define the highest form of a 
nilpotent matrix and then  
reduce the general case to a sum of nilpotents.  We begin with some
 linear algebra.

\begin{definition}
Fix a matrix $X$.  The entry $X_{ik}$ is a pivot of $X$ if $X_{ik}$ is 
nonzero and if all entries below and to its left vanish, that is
$X_{ij}=0$ if $j < k$ and $X_{jk}=0$ if $j > i$.
\end{definition}

Given $i$, define $r_i$ to be the row of $X_{r_i,i}$ 
if the entry is a pivot and zero if not.

\begin{definition}
Fix an upper-triangular nilpotent matrix $N$.  Then $N$ is in 
highest form if the pivots form a nondecreasing sequence, namely
$r_1 \leq r_2 \leq \cdots \leq r_n$.
\end{definition}

By definition $r_i = r_j$ only if both are zero, so only initial columns
of a matrix in highest form can be zero.  
Columns with pivots are linearly independent, so when $N$ is in highest 
form its first $\dim(\ker N)$ columns are zero.

To construct a highest form for $N$ fill the Young
diagram $\lambda_N$ constructed in the Introduction with $1$ to $n$ starting
at the bottom of the leftmost column, incrementing by one while moving up, 
then moving to the lowest box of the next column and repeating.  
The highest form for $N$ is the matrix with $N_{ij}=1$
if $i$ fills the box to the left of $j$ and $N_{ij}=0$ otherwise, as
in Figure \ref{highest form example}.

\begin{figure}[h]
\[\mbox{\small $\left( \begin{array}{ccccccc}
0 & 1 & 0 & 0 & 0 & 0\\
0 & 0 & 1 & 0 & 0 & 0\\
0 & 0 & 0 & 0 & 0 & 0\\
0 & 0 & 0 & 0 & 1 & 0\\
0 & 0 & 0 & 0 & 0 & 0\\
0 & 0 & 0 & 0 & 0 & 0
\end{array} \right) $}
\longleftrightarrow \hspace{.2in}
\textup{\begin{tabular}{|c|c|c|}
\cline{1-3} 3 & 5 & 6 \\
\cline{1-3} 2 & 4 & \multicolumn{1}{c}{} \\
\cline{1-2} 1 & \multicolumn{2}{c}{} \\
\cline{1-1} \multicolumn{3}{c}{} \end{tabular}} 
\longleftrightarrow
\mbox{\small$\left( \begin{array}{ccccccc}
0 & 0 & 0 & 0 & 0 & 0\\
0 & 0 & 0 & 1 & 0 & 0\\
0 & 0 & 0 & 0 & 1 & 0\\
0 & 0 & 0 & 0 & 0 & 0\\
0 & 0 & 0 & 0 & 0 & 1\\
0 & 0 & 0 & 0 & 0 & 0
\end{array} \right)$}\]
\caption{Jordan canonical form, the Young diagram,
and the highest form} \label{highest form example}
\end{figure}

The main property of the highest form is that conjugation by $U$
preserves it.

\begin{proposition} \label{nilpotent highest form preserved}
If $N$ is nilpotent and 
in highest form and $u \in U$ then $u^{-1}Nu$ is in highest
form.  The entry $N_{r_j,j}$ is a pivot if and only if 
$(u^{-1}Nu)_{r_j,j}$ is.  If so, $N_{r_j,j}=
(u^{-1}Nu)_{r_j,j}$.
\end{proposition}

\begin{proof}
The entry $(Nu)_{jk}$ is the sum of $N_{jk}$ and 
multiples of $N_{j1}$, $\ldots$, $N_{jk-1}$.  This means $(Nu)_{jk}=
N_{jk}$ for each column up to and including the first
nonzero column in the $j^{th}$ row of $N$.  Similarly 
$(u^{-1}Nu)_{jk}=(Nu)_{jk}$ for each row
after and including the last nonzero row in the $k^{th}$ column of $Nu$.
So the pivots of $u^{-1}Nu$ are in the same entries with
 the same values as in $Nu$, which are in the same entries
with the same values as in $N$.
\end{proof}

We now describe highest form for an arbitrary upper-triangular matrix
$S+N$, where $S$ is diagonal and $N$ is nilpotent.
If $c$ is an eigenvalue of $S$ then let $E_c$ be its eigenspace.
Recall that $S$ induces a decomposition of the total vector space
$\C^n=\bigoplus_{\textup{eigenvalues $c$ of $S$}} E_c$.

Inclusion and then projection gives a map from the semigroup
$\End(\C^n)$ to $\End(E_c)$.  For instance, the image of $S+N$ under
this map is the composition
\[E_c \hookrightarrow \C^n \stackrel{S+N}{\longrightarrow} 
   \C^n \longrightarrow 
   {\left( \C^n / \bigoplus_{c' \neq c} E_{c'} \right)} 
   \hspace{.05in} \cong E_c.\]
The matrix for $(S+N)_c$ is given by 
the $\dim E_c \times \dim E_c$
minor of $S+N$ obtained by removing the $j^{th}$ row and $j^{th}$ column
if $S_{jj} \neq c$.  This is shown in Figure \ref{submatrix example}.

\begin{figure}[h]
\[S+N = \left( \begin{array}{ccc} 1 & a & b \\ 0 & 0 & c \\ 0 & 0 & 1
    \end{array} \right) \hspace{.2in} \mapsto \hspace{.2in} (S+N)_1 = 
    \left( \begin{array}{cc} 1 & b \\ 0 & 1 \end{array} \right)\]    
\caption{An example of $S+N$ and $(S+N)_1$}\label{submatrix example}
\end{figure}

%
%
Note that $N_c$ is the strictly upper-triangular part of $(S+N)_c$.

\begin{definition}
$S+N$ is in highest form if the following hold:
\begin{enumerate}
\item $S+N$ is upper triangular;
\item \label{diag condition} 
  if $S_{ii}=S_{jj}$ then $S_{ii}=S_{kk}$ for each $k$ between $i$ and $j$;
  and
\item $N_c$ is in highest form for each eigenvalue $c$ of $S$.
\end{enumerate}
\end{definition}

The diagonal blocks of a matrix in highest form are in
highest form.  However, highest form matrices 
need not be block diagonal in general.
Condition \ref{diag condition} is designed so the map
$Y \mapsto Y_c$ is a morphism of semigroups, as in 
the next lemma.  Again $e_i$ is the standard basis vector in $\C^n$.

\begin{lemma} \label{morphism condition}
$(XY)_c=X_cY_c$ for all upper-triangular matrices $X$ and $Y$ if and only
if there are $i$ and $j$ so that
$E_c$ is the span of the basis vectors $e_i$,$e_{i+1}$,$e_{i+2}$,$\ldots$,
$e_{i+j}$.
\end{lemma}

\begin{proof}
The coefficient of $e_k$ in $XY e_{i'}$
is $\sum_{j=k}^{i'} x_{kj} y_{ji'}$.
If $E_c$ satisfies the hypothesis then for each $e_k$ spanning $E_c$
the entries $x_{kj}$ and
$y_{ji'}$ are in $X_c$ and $Y_c$ respectively as long as $j$ is between
$k$ and $i'$.  Consequently $(XY)_c=X_cY_c$.

Conversely, suppose $e_i$,
$e_{i+k}$, and $e_{i+j}$ are vectors 
with $0<k<j$ and with $e_i, e_{i+j}$ in $E_c$
while $e_{i+k}$ is not. If $X$ is a matrix nonzero only in entry
 $X_{i,i+j}$ and $Y$ is nonzero only in entry
 $Y_{i+j,i+k}$ then 
$X_c=Y_c=0$ but $(XY)_c$
is nonzero.
\end{proof}

To construct a matrix in highest form which is conjugate to $S+N$,
 write $S+N$ in Jordan canonical form
$\sum (S_i + N_i)$ with blocks $S_i+N_i$ 
corresponding to distinct eigenvalues $c_i$.  
If $N_i'$ is highest form for $N_i$ then the matrix 
$\sum (S_i + N_i')$ is in highest form, called the
permuted Jordan form of $S+N$.

The next proof extends Proposition \ref{nilpotent highest form preserved} 
to general linear operators.

\begin{proposition} \label{semisimple highest form preserved}
If $S+N$ is in highest form and $u$ is in $U$ then $u^{-1}(S+N)u$ 
is in highest form.   
The $(r_k,k)$ entry of $N_c$ is a pivot if and only if the
$(r_k,k)$ entry of $(u^{-1}Nu)_c$ is a pivot.  
In this case the two entries are equal.
\end{proposition}

\begin{proof}
Note that $u^{-1}(S+N)u$ is upper-triangular if $S+N$ is.  
Direct computation shows $u^{-1}(S+N)u = S+(u^{-1}Nu + u^{-1}Su - S)
= S+N'$ for some nilpotent $N'$. 

Fix an eigenvalue $c$ of $S$.
By Lemma \ref{morphism condition} 
we know $(u^{-1}Nu)_c = u^{-1}_cN_c u_c$.  
Proposition \ref{nilpotent highest form preserved} applies
since $N_c$ is in highest form and $u_c$ is 
upper-triangular with ones on the diagonal.
\end{proof}

\section{Paving Hessenberg varieties by affines} \label{paving by affines section}

In this section we prove that if
$X$ is in highest form, each row of each Schubert cell is in $\H(X,h)$
if and only if certain affine conditions hold.  This is the key
step in the paper.

Recall that $\vect{X}{i}$ is the $i^{th}$ row of $X$, that 
$\colvect{X}{j}$ is the $j^{th}$ column, and $H$ is the Hessenberg
space given by $h$ in Definition \ref{Hessenberg space}.
The next lemma identifies $\{u \in U_i : \vect{\left(u^{-1}Nu\right)}{i}
\in \vect{\left(wHw^{-1}\right)}{i}\} \cap U_w$ 
as the solution to an affine system of equations and finds its rank.

\begin{lemma} \label{nilpotent row}
Fix a permutation $w$, a row $U_i$, 
a Hessenberg space $H$, and $N$ in highest form.
If the pivots of $N$ are in nonzero entries of $wHw^{-1}$ 
then the set $\{u \in U_i: \vect{\left(u^{-1}Nu\right)}{i}
\in \vect{\left(wHw^{-1}\right)}{i}\} \cap U_w$ 
is homeomorphic to $\C^d$ for 
\[d=|\{k: k>i,  w(i) > w(k), h(w(j)) \geq w(i) \textup{ if } 
   N_{kj} \textup{ is a pivot in } N\}|.\]
\end{lemma}

The inequality $h(w(j)) \geq w(i)$ does not apply
if the $k^{th}$ row of $N$ has no pivot.

\begin{proof}
The $i^{th}$ row of $u^{-1}Nu$
is $\vect{(u^{-1})}{i}N$ by Proposition \ref{affine transformation}. 
Examining the 
condition $\vect{(u^{-1})}{i}N \in \vect{\left(wHw^{-1}\right)}{i}$
for each column gives the system of equations
\[ \vect{(u^{-1})}{i}\colvect{N}{j}=0 \hspace{.25in} \textup{ for } j
       \textup{ such that }
          w(i)>h(w(j)).  \]
Each equation in this system has the form 
\[(1, -u_{i,i+1}, \ldots, -u_{i,n}) \cdot (N_{i,j}, \ldots, 
N_{n,j})^t = 0\]
for $j$ satisfying $w(i)>h(w(j))$.  Adding the constraint that
$u \in U_w$ gives
the following affine system of equations in the free entries $u_{ik}$:
\begin{equation} \label{affine system}
  (u_{ik_1}, u_{ik_2}, \ldots, u_{ik_{d_i}}) 
   \left( \begin{array}{c} N_{k_1j} \\
     N_{k_2j} \\ \vdots \\ N_{k_{d_i}j} \end{array} \right) = N_{ij} 
    \hspace{.2in} \begin{array}{c}
      \textup{ for } j \textup{ with } w(i)>h(w(j)) \\
    \textup{ and } k_l \textup{ with }w(i) > w (k_l).\end{array}
  \end{equation}
The linear system of equations ${\bf x} M = {\bf v}$ 
has a solution if and only if 
the rank of the coefficient matrix $M$ equals that of the 
extended matrix $\binom{{\bf v}}{M}$.
To prove this here, we show that if either $N_{ij}$ or
one of the $N_{k_lj}$ is nonzero
then in fact one of the $N_{k_lj}$ is a pivot in $N$.

Indeed, if $N_{ij}$ or $N_{k_lj}$ is nonzero
 then $N$ has a pivot $N_{kj}$ in some row
$k \geq i$.  The pivots of $N$ are in $wHw^{-1}$ by hypothesis.
This means that $w(k) \leq h(w(j))$ by Proposition \ref{algebraic to 
combinatorial}.  In addition $w(i)>h(w(j))$ by hypothesis on $j$.  
Hence $w(i) > w(k)$ and 
so $N_{kj}$ is one of the entries of the column vector of 
Equation \eqref{affine system}.

The dimension of the solution space is the 
number of free entries in $U_i \cap U_w$
less the number of pivots of $N$ in the coefficient matrix of Equation
\eqref{affine system}.  The set $\{k: k > i, w(i) > w(k)\}$ indexes
the free entries while $\{k: k>i, w(i) > w(k), 
  N_{kj} \textup{ is a pivot and } w(i)>h(w(j))\}$
indexes the rank of the coefficient matrix.  This proves the claim. 
\end{proof}

This extends to general linear operators in much the same way.

\begin{lemma} \label{semisimple row}
Fix a permutation $w$, a row $U_i$, 
a Hessenberg space $H$, and $S+N$ in highest form.
If the pivots of each submatrix $N_c$ are in $wHw^{-1}$ 
then the set $\{u \in U_i: \vect{\left(u^{-1}(S+N)u\right)}{i}
\in \vect{\left(wHw^{-1}\right)}{i}\} \cap U_w$ 
is homeomorphic to $\C^d$ for 
\begin{eqnarray*} d & = &|\{k: k>i,w(i) > w(k),  \\
   && \hspace{.3in}h(w(j)) \geq w(i) 
  \textup{ if } 
  N_{kj} \textup{ is a pivot in } N_{S_{ii}}, 
   S_{kk}=S_{ii}\}| \\ 
  & + & |\{k: k>i, h(w(k)) \geq w(i) > w(k), 
   S_{kk} \neq S_{ii}\}|. 
\end{eqnarray*}
\end{lemma}

\begin{proof}
The $i^{th}$ row of $u^{-1}(S+N)u$ is 
$S_{ii}\vect{u}{i} + \vect{(u^{-1})}{i}(S+N)$ by Proposition 
\ref{affine transformation}. 
The condition that this be in $\vect{\left(wHw^{-1}\right)}{i}$
 gives the system of equations
\[ S_{ii}u_{ij}+ \vect{(u^{-1})}{i}\colvect{(S+N)}{j}=0 
      \hspace{.25in} \textup{ for } j
       \textup{ such that }
          w(i)>h(w(j)).  \]
Each equation in this system is of the form 
\[S_{ii}u_{ij} + (1, -u_{i,i+1}, \cdots, -u_{i,n}) \cdot (N_{ij},
\cdots, N_{j-1,j}, S_{jj}, 0, \ldots, 0)^t = 0\] 
for $j$ such that $w(i)>h(w(j))$.  Adding the condition that $u \in U_w$  
gives the system
\[
  (u_{ik_1}, u_{ik_2}, \ldots, u_{ik_{d_i}}) 
   \left( \begin{array}{c} N_{k_1j} \\
     N_{k_2j} \\ \vdots \\ 
     S_{jj}-S_{ii} \\ 0 \\ \vdots \\ 0
   \end{array} \right) = 
   N_{ij} 
    \hspace{.25in} \begin{array}{c}
      \textup{ for } j \textup{ such that } w(i)>h(w(j)) \\
    \textup{ and } k_l \textup{ such that }w(i) > w (k_l).\end{array}
  \]
As in the previous lemma, we show that the rank of the coefficient
matrix is unchanged if the vector of solutions $(N_{ij})$ is inserted
as the top row.

We study the cases when $S_{ii}=S_{jj}$ and when $S_{ii} \neq S_{jj}$
separately.  Let $c_i$ be the cardinality $|\{S_{jj}: j>i, S_{jj}=
S_{ii}\}|$ so $S_{jj}-S_{ii}$
is zero exactly when $j$ is at most $i+c_i$.  
The columns with $j > i+c_i$ have a pivot in
position $(j,j)$ regardless of $N_{ij}$.
For each such $j$ we know $w(i) > w(j)$ since $h(w(j)) \geq w(j)$.

The first $c_i$ columns and rows of this system satisfy $S_{jj}-S_{ii}=0$
and so form the system of Equation \eqref{affine system}.
Its pivots are computed in Lemma \ref{nilpotent row}.  
Each is a pivot in the original system because the $(k_l,j)^{th}$
entry is zero when $k_l$ is greater than $j$.

The rank of the entire matrix is therefore
\[\begin{array}{l}|\{k: k>i, w(i) > w(k), w(i) > h(w(j)) \textup{ and } 
N_{kj} \textup{ is a pivot in } N_{S_{ii}}, S_{ii}=S_{jj}\}| \\
\textup{  } +
|\{k: k>i,w(i) > w(k), w(i) > h(w(k)), S_{ii} \neq S_{kk}\}|.\end{array}\]
  Since 
the dimension of $U_i \cap U_w$ is $|\{k: k>i, w(i)>w(k)\}|$ the 
claim follows.
\end{proof}

\section{The Main Theorems} \label{main theorems}

We now demonstrate that requiring each
row of a flag in $\H(X,h)$ to satisfy the Hessenberg conditions 
gives the structure of 
an iterated tower of affine fiber bundles on each Bruhat cell in
$\H(X,h)$.  This constructs a paving by affines on the Hessenberg variety.
We use the Hessenberg space $H$ determined by $h$ as in Definition 
\ref{Hessenberg space}, as well as the description of the Schubert cells
in Proposition \ref{Bruhat decomposition}.

\begin{theorem} \label{main theorem: bundle form}
Fix a Hessenberg space $H$ and a basis for which $S+N$ is 
in highest form and  in permuted Jordan form.  Let $\{C_w\}$ be the Schubert
cells.

The intersections $C_w \cap \H(S+N,H)$ form a paving by affines of
$\H(S+N,H)$.
The cell $C_w \cap \H(S+N,H)$ is nonempty 
if and only if $N$ is in $wHw^{-1}$.  If nonempty, the cell
 $C_w \cap \H(S+N,H)$ is homeomorphic to $\mathbb{C}^d$
for
\begin{eqnarray*} d &=& |\{(i,k): k>i, w(i) > w(k), \\
  & & \hspace{.5in} h(w(j)) \geq w(i) \textup{ if }
  N_{kj} \textup{ is a pivot in } N_{S_{ii}}, S_{kk}=S_{ii}\}| 
  \\ &+& |\{(i,k): k>i, h(w(k)) \geq w(i) > w(k), 
   S_{kk} \neq S_{ii}\}|. \end{eqnarray*}
\end{theorem}

\begin{proof}
The Schubert cells $\{C_w\}$ form a paving of the full flag variety.
The Hessenberg variety $\H(S+N,H)$ is a closed subvariety of the 
flag variety so the intersections $C_w \cap \H(S+N,H)$ pave the
Hessenberg variety.

We now identify the nonempty cells.
If $N$ is in $wHw^{-1}$ then the flag $[w]$
is in $\H(S+N,H)$.  Conversely, if
the flag $[uw]$ is in $\H(S+N,H)$ then $w^{-1}
u^{-1}(S+N) uw \in H$.  This implies that the pivots of each submatrix 
$(u^{-1}Nu)_{S_{ii}}$ are in $w H w^{-1}$.  Since $N$ is in highest form,
its pivots are in the same positions as those of $u^{-1}Nu$
 by Proposition \ref{nilpotent highest form preserved}.  
Each pivot of $N$ is a pivot of some $N_{S_{ii}}$ because
$S+N$ is in permuted Jordan form. 
The pivots of each $N_{S_{ii}}$ are in $w H w^{-1}$ if and only if those of 
$(u^{-1}Nu)_{S_{ii}}$ are.  The only nonzero entries of $N$ are 
pivots so $N$ is in $wHw^{-1}$.

Next, suppose $C_w \cap \H(S+N,H)$ is nonempty.  
Define 
\[\begin{array}{l} Z_{i} = \left\{ u \in \left(U_{n-1}^{} U_{n-2} 
      \cdots U_{i}\right) \vspace{1.5em} \cap U_w : \right. \\
\hspace{.75in} \left. \vect{\left({u}^{-1}(S+N)u\right)}{j}
\in \vect{\left(wHw^{-1}\right)}{j} \textup{ for all } j > i \right\}.
\end{array}\]  
For instance, $Z_{n-1} = U_{n-1} \cap U_w$ since $wH w^{-1}$ always contains
the span of $E_{nn}$.  Also, observe that $Z_1$
is homeomorphic to $C_w \cap \H(S+N,H)$ under the map which sends
$u \mapsto uw$.  
We will show that $Z_1$ is affine and compute its dimension.

To do this, we factor each element in $Z_i$ uniquely as $u'u$ for
$u' \in U_{n-1} \cdots U_{i+1}$ and $u \in U_i$ by Proposition 
\ref{rows factor U}.  Conjugation by $U_i$ only
affects the first $i$ rows of an upper triangular matrix by 
Proposition \ref{row conjugation}, so $u^{-1}({u'}^{-1}(S+N)u')u$ agrees
with ${u'}^{-1}(S+N)u'$ in rows $i+1$ and higher.  
Thus, this factorization satisfies the additional
conditions that $u' \in Z_{i+1}$ and that $u \in U_i \cap U_w$
has $\vect{\left(u^{-1}({u'}^{-1}(S+N)u')u\right)}{i}
\in \vect{\left(wHw^{-1}\right)}{i}$.  This gives a well-defined
map $\pi_i: Z_i \rightarrow Z_{i+1}$ sending $u'u$ to $u'$.

We now show that $\pi_i: Z_i \rightarrow Z_{i+1}$ 
is an affine fiber bundle and
compute its rank.  For each element $u' \in Z_{i+1}$, the operator
${u'}^{-1} (S+N) u'$ is in highest form and has its pivots in the same
position as $S+N$.  
Consequently, the hypotheses of Lemma \ref{semisimple row}
hold.  Lemma \ref{semisimple row} states that for each $u' \in 
Z_{i+1}$, the preimage $\pi_i^{-1}(u') \subseteq Z_i$ 
is affine of dimension
\begin{eqnarray*} d_i &=&|\{k: k>i,  w(i) > w(k), \\
  && \hspace{.3in}h(w(j)) \geq w(i) 
  \textup{ if } N_{kj} \textup{ is a pivot in } N_{S_{ii}}, 
   S_{kk}=S_{ii}\}| \\ & + &
   |\{k: k>i, h(w(k)) \geq w(i) > w(k), 
   S_{kk} \neq S_{ii}\}|. \end{eqnarray*} 
The fiber $\pi_i^{-1}(u')$ is the set of solutions ${\bf x}_{u'}$ to the
affine system ${\bf x}_{u'} M_{u'} = {\bf v}_{u'}$, where $M_{u'}$ 
and ${\bf v}_{u'}$ vary continuously (by conjugation) in $u'$.  In other
words $\pi_i: Z_i \longrightarrow Z_{i+1}$ is a fiber bundle.

We produce a bundle homeomorphism from
$\pi_i: Z_i \longrightarrow Z_{i+1}$ to the trivial bundle of rank $d_i$
over $Z_{i+1}$.  Let $I$ be the set of indices used to define $d_i$ in 
Lemma \ref{semisimple row}.  For each $u' \in Z_{i+1}$, Lemma
\ref{semisimple row} shows that the $(i,k)$ entry of the matrices
in $\pi_i^{-1}(u')$ is free whenever $k \in I$.
The map sending $u'u \mapsto (u', (u_{ik})_{k \in I})$
has a continuous inverse given by the system ${\bf x}_{u'} M_{u'} = 
{\bf v}_{u'}$ and so is the desired bundle homeomorphism.  
Given this bundle map, if
$Z_{i+1}$ is homeomorphic to affine space then $Z_i$ is homeomorphic to
affine space of dimension $\dim Z_{i+1} + d_i$.

Finally, consider the sequence $Z_1 \stackrel{\pi_1}{\longrightarrow} Z_2
\stackrel{\pi_2}{\longrightarrow} Z_3 \cdots 
\stackrel{\pi_{n-2}}{\longrightarrow} Z_{n-1}$.  Each map $\pi_i$ is
an affine fiber bundle of rank $d_i$.  We know $Z_{n-1} = U_{n-1} \cap U_w$
is affine and write 
its dimension as $d_{n-1} = |\{k: k> n-1, w(n-1) > w(k)\}|$ to stress
the analogy to the other $d_i$.  
Inducting on $i$, we may assume the base space of 
$Z_i \stackrel{\pi_i}{\longrightarrow} Z_{i+1}$ 
is homeomorphic to affine space, and so
its total space $Z_i$ is homeomorphic to affine space of dimension
$\dim Z_{i+1} + d_i$.
By induction $Z_1$ is homeomorphic to $\mathbb{C}^d$ with $d=d_1 + 
\cdots + d_{n-1}$.  
\end{proof}

This along with Proposition \ref{betti of pavings} leads 
to an immediate corollary when the base field is $\C$.

\begin{corollary}
Hessenberg varieties have no odd-dimensional cohomology.
\end{corollary}

The main theorem is much simpler if the operator is 
nilpotent or semisimple.  

\begin{corollary} \label{nilpotent: bundle form}
Fix a Hessenberg space $H$.  Let $N$ be a nilpotent matrix in highest form
and in permuted Jordan form.  Let $\{C_w\}$ be the Schubert
cells.

The intersections $C_w \cap \H(N,H)$ form a paving by affines of
$\H(N,H)$.
The cell $C_w \cap \H(N,H)$ is nonempty 
if and only if $N$ is in $wHw^{-1}$.  If nonempty,
the cell $C_w \cap \H(N,H)$ is homeomorphic to $\mathbb{C}^d$
for
\[d = |\{(i,k): k>i, w(i) > w(k), h(w(j)) \geq w(i) \textup{ if }
  N_{kj} \textup{ is nonzero}\}|.\]
\end{corollary}

The proof of this is immediate, as is that of the next corollary.

\begin{corollary}
Fix a Hessenberg space $H$.  Let $S$ be a diagonal matrix in highest form
and let $\{C_w\}$ be the Schubert cells of the flag variety.  
The intersections $C_w \cap \H(S,H)$ form a paving by affines of
$\H(S,H)$.  The cell $C_w \cap \H(S,H)$
is homeomorphic to $\mathbb{C}^d$ for
\[d = |\{(i,k): k>i, w(i) > w(k),  h(w(k)) \geq w(i) 
   \textup{ if } S_{kk} \neq S_{ii}\}|.\]
\end{corollary}

In particular, the intersection of 
each Schubert cell with $\H(S,H)$
is nonempty!  

\begin{corollary}
If $S$ is diagonal then the Euler characteristic $\chi(\H(S,h))$
is $n!$ for every Hessenberg function $h$.
\end{corollary}

\begin{proof}
Since $w^{-1}Sw$ is diagonal for each permutation, every Schubert cell
$C_w$ intersects $\H(S,h)$ in a nonempty affine cell 
$\C^{d_w}$.  Since the cohomology
is only even-dimensional, the Euler characteristic of $\H(S,h)$ is
the total number of cells.
\end{proof}

\section{Tableaux Interpretations} \label{main theorem: tableau form}

We describe the main theorems combinatorially using Young
diagrams.

To each linear operator $X$ we associate a multitableau $\lambda_X$
as follows.
If $\sum (S_i + N_i)$ is a Jordan canonical form for $X$
then $\lambda_X$ is the collection of tableaux $\lambda_{N_i}$ 
associated to $N_i$ as in the Introduction.
We assume tableaux are ordered vertically by size 
as shown in Figure \ref{defn multitableaux}.  
Note that $\lambda_X$ is independent of the numerical eigenvalues 
of $S_i$.  When $X$ is nilpotent this definition
reduces to that of Figure \ref{defn tableaux}.

\begin{figure}[h]
\[\mbox{\small{$\left( \begin{array}{ccccccc}
0 & 1 & 0 & 0 & 0& 0& 0\\
0 & 0 & 1 & 0 & 0& 0& 0\\
0 & 0 & 0 & 0 & 0& 0& 0\\
0 & 0 & 0 & 0 & 0& 0&0 \\
0 & 0 & 0 & 0 & 2& 1&0 \\
0 & 0 & 0 & 0 & 0& 2&0 \\
0 & 0 & 0 & 0 & 0& 0&2 
\end{array} \right)$}}
\hspace{.05in} \longleftrightarrow \hspace{.05in}
\textup{\begin{tabular}{l}
{\begin{tabular}{|c|c|c|}
\cline{1-3} 5 & 6 & 7 \\
\cline{1-3} 4 & \multicolumn{2}{c}{} \\
\cline{1-1} \multicolumn{3}{c}{}  \end{tabular}} \\
{\begin{tabular}{|c|c|}
\cline{1-2} 2 & 3 \\
\cline{1-2} 1 & \multicolumn{1}{c}{} \\
\cline{1-1} \multicolumn{2}{c}{} \end{tabular} } \end{tabular}} 
\hspace{.05in} \longleftrightarrow \hspace{.05in}
\mbox{\small{$\left( \begin{array}{ccccccc}
 2& 0&0 & 0 & 0 & 0 & 0 \\
 0& 2&1 & 0 & 0 & 0 & 0 \\
 0& 0&2 & 0 & 0 & 0 & 0 \\
0 & 0 & 0 & 0 & 0 & 0 & 0  \\
0 & 0 & 0 & 0 & 0 & 1 & 0  \\
0 & 0 & 0 & 0 & 0 & 0 & 1  \\
0 & 0 & 0 & 0 & 0 & 0 & 0  
\end{array} \right)$}}\]
\caption{The Jordan form, multitableau with base filling, and highest
form of a general linear operator}\label{defn multitableaux}
\end{figure}

The base filling of $\lambda_X$
is that for which each $\lambda_{N_i}$ is filled according to the rules  
in Figure \ref{highest form 
example} except that the lowest number in 
$\lambda_{N_i}$ is one more than 
the highest in $\lambda_{N_{i-1}}$.
Figure \ref{defn multitableaux} demonstrates
this.  The box containing $i$ in this filling
of $\lambda_X$ is called the $i^{th}$ box.

We associate each filling of the multitableau $\lambda_X$ to 
a unique permutation $w$ according to the convention that the $i^{th}$ box
contains $w(i)$.  For instance, the $i^{th}$ box of the base filling
contains $i$.

\begin{theorem} \label{theorem: tableau form}
Fix any linear operator $X$ and Hessenberg function $h$.  
The Hessenberg variety 
$\H(X,h)$ is paved by affines.  The nonempty
cells are naturally in bijection with the
fillings of $\lambda_X$ which contain the configuration 
\begin{tabular}{|c|c|} \cline{1-2} $k$ & $j$ \\ \cline{1-2} 
\end{tabular} only if $k \leq h(j)$.  The dimension of a nonempty
cell is the sum of:

\begin{enumerate}
\item \label{nilpotent condition} the number of pairs $i$,$k$
in the corresponding filling of $\lambda_X$ such that 
\begin{itemize}
\item $i$ and $k$ are in the same tableau,
\item the box filled by $i$ is to the left of or directly below 
  the box filled by $k$, 
\item $k < i$, and
\item  if $j$ fills the box immediately to the right of 
  $k$ then $i \leq h(j)$.
\end{itemize}
\item \label{semisimple condition} 
the number of pairs $i$,$k$ in $\lambda_X$ such that
\begin{itemize}
\item $i$ and $k$ are in different tableaux,
\item the box filled with $i$ is below $k$, and
\item $k < i \leq h(k)$.
\end{itemize}
\end{enumerate}
\end{theorem}

The first condition is illustrated
in Figure \ref{picture of conditions} 
and the second in Corollary 
\ref{semisimple tableau}.

\begin{proof}
Write $i'$ for the index of the box containing $i$, respectively
$j'$ and $k'$.  This means that $w(i') = i$ so $i>k$ if 
and only if $w(i') > w(k')$.  

The $i'^{th}$ box is in the same tableau
as the $k'^{th}$ box if and only if $S_{i'i'}=S_{k'k'}$.

Box $i'$ sits left of or directly below box $k'$ if and only if $k' > i'$
by the labelling convention.  

The nilpotent part of a permuted Jordan form is the sum of $E_{k'j'}$
over $(k',j')$ such that box $j'$ sits to the right of box $k'$.  $X$ is in
$wHw^{-1}$ exactly when each of these summands is and
each $E_{k'j'}$ is in $wHw^{-1}$ exactly when $k=w(k') \leq h(w(j'))=h(j)$
by Proposition \ref{algebraic to combinatorial}.
\end{proof}

We prove Theorem \ref{nilpotent: tableau form}, paving 
nilpotent Hessenberg varieties using tableaux.

\begin{proof}
If $N$ is nilpotent its  multitableau consists of
exactly one tableau.  Condition \ref{semisimple condition} of Theorem
\ref{theorem: tableau form} never applies so 
Condition \ref{nilpotent condition} gives the dimension.
\end{proof}

The following interprets the main theorem for semisimple operators.

\begin{corollary} \label{semisimple tableau}
Fix a Hessenberg space $H$.  Let $S$ be a diagonal matrix and $\lambda_S$
its associated multitableau.  
The Schubert cell $C_w$ intersects the
Hessenberg variety $\H(S,H)$ in a space homeomorphic to $\mathbb{C}^d$
where $d$ is the sum of:
\begin{enumerate}
\item the number of pairs $i$,$k$ such that  

  \vspace{-.25in} 
\hspace{2in} (1) \begin{tabular}{|c|}
    \multicolumn{1}{c}{$\vdots$} \\
    \cline{1-1} $\vdots$ \\
    \cline{1-1} $k$ \\
    \cline{1-1} $\vdots$ \\
    \cline{1-1} $i$ \\
    \cline{1-1} $\vdots$ \\
    \cline{1-1} \multicolumn{1}{c}{$\vdots$}
    \end{tabular} 

\vspace{-1.35in}
\begin{itemize}
\item $i$ and $k$ are in the same tableau,
\item $i$ is below $k$, and
\item $k < i$.
\end{itemize}
\item the number of pairs $i$,$k$ such that 

\vspace{-0.95in} 
\hspace{3in}  (2) \begin{tabular}{|c|}
    \cline{1-1} $\vdots$ \\
    \cline{1-1} $k$ \\
    \cline{1-1} $\vdots$ \\
    \cline{1-1} \multicolumn{1}{c}{\small $\vdots$} \\
    \cline{1-1} $\vdots$ \\
    \cline{1-1} $i$ \\
    \cline{1-1} $\vdots$ \\
    \hline
    \end{tabular} 

\vspace{-.65in}
\begin{itemize}
\item $i$ and $k$ are in different tableaux,
\item $i$ is below $k$, and
\item $k < i \leq h(k)$.
\end{itemize}
\end{enumerate}
\end{corollary}

\vspace{.05in}
\begin{proof}
The nilpotent associated to each eigenspace is
the zero matrix so each Young diagram is a single column.  
This implies that every Schubert cell intersects the 
Hessenberg variety and that the first
condition of Theorem \ref{theorem: tableau form} simplifies as given.
\end{proof}

\section{Root system interpretation} \label{main theorem: root form}

The main theorem can also be expressed in terms of roots.  For general
background on Lie algebras, the reader is referred to \cite{H2}.  

Recall that the Lie algebra of $GL_n(\C)$
is $\mathfrak{gl}_n(\C)$, which we think of as $n \times n$ matrices
over $\C$.  Fix the Borel
subalgebra $\mathfrak{b}$ of upper-triangular matrices
in $\mathfrak{gl}_n(\C)$.  

The standard embedding of $\mathfrak{gl}_n(\C)$
into the space of matrices associates the matrix $E_{ij}$ with $i < j$ 
to the root vector $E_{\alpha}$ where
$\alpha = \alpha_i + \alpha_{i+1} + \ldots + \alpha_{j-1}$.  
The root $\alpha$ can also be regarded as the linear functional 
on diagonal matrices with $\alpha(S)=S_{jj}-S_{ii}$.

The set of positive roots $\Phi^+$ are the roots $\alpha$
for which $E_{\alpha}$ is upper-triangular.  
The set of negative roots $\Phi^-$ are the roots $-\alpha$
for $\alpha$ in $\Phi^+$.  They correspond to the lower-triangular 
matrices by the map which sends $E_{ji}$ to $-\alpha$
if $E_{\alpha}=E_{ij}$.  The action of the permutation $w$ on the set of
roots is defined by $w^{-1} \alpha = \beta$ 
if $w^{-1}E_{\alpha}w = E_{\beta}$.

With this notation
a Hessenberg space $H$ can be defined intrinsically as a vector subspace
of $\mathfrak{gl}_n(\C)$ which contains $\mathfrak{b}$ and which is
closed under Lie bracket with $\mathfrak{b}$ as in \cite{dMPS}.  
We write $\Phi_H$ to
denote the roots whose root spaces span $H$.

The definition of highest form operators can be extended to root spaces
by the standard embedding.
If $S+N$ is in highest form we denote by $\Phi_{S+N}$ the set of roots
corresponding to the pivots of $N_c$ over all eigenvalues $c$ of $S$.

\begin{theorem} \label{theorem: root form}
Fix a Hessenberg space $H$.  Fix $\mathfrak{b}$ 
with respect to which $S+N$ is 
in highest form and permuted Jordan form.  
The intersection $C_w \cap \H(S+N,H)$ is nonempty if and 
only if $w^{-1} \Phi_{S+N}$ is in $\Phi_H$.  
If so $C_w \cap \H(S+N,H)$ is 
homeomorphic to $\mathbb{C}^d$ for
\begin{eqnarray*} d &=& |\{\alpha \in \Phi^+: 
\alpha(S)=0, w^{-1}\alpha \in \Phi^-,
  w^{-1}(\alpha+\beta) \in \Phi_H \textup{ for some } \beta \in \Phi_{S+N}\}| 
  \\ &+& |\{\alpha \in \Phi^+: \alpha(S) \neq 0, w^{-1} \alpha \in \Phi_H,
  w^{-1} \alpha \in \Phi^-\}|. \end{eqnarray*}
\end{theorem}

\begin{proof}
Write $N$ in terms of root vectors as 
$\sum_{\beta \in \Phi_{S+N}} E_{\beta}$.

The pivot $E_{\beta}$ is in $wHw^{-1}$
if and only if $w^{-1} \beta \in \Phi_H$ by Proposition 
\ref{algebraic to combinatorial}.

If $\alpha = \alpha_i + \alpha_{i+1} + \cdots + \alpha_{k-1}$
then $S_{ii}=S_{kk}$ if and only if $\alpha(S)=0$, which describes
two of the conditions in the theorem.  

The condition $h(w(k)) \geq w(i)$ is
equivalent to $w^{-1} 
\alpha \in \Phi_H$  by Proposition \ref{algebraic to
combinatorial}.

The root $\alpha$ satisfies
$k>i$ and $w(i)>w(k)$ if and only if $\alpha \in \Phi^+$ and 
$w^{-1} \alpha \in \Phi^-$ according to 
the characterization of the Bruhat decomposition in Proposition
\ref{Bruhat decomposition}.  

The condition that $N_{kj}$ be a pivot in $N_{S_{ii}}$ 
indicates that $\beta = \alpha_k + \alpha_{k+1} + \cdots + \alpha_{j-1}$
is a root in $\Phi_{S+N}$.  The root $\alpha + \beta$ corresponds to 
$E_{ij}$.  This means that 
the condition $w^{-1}(\alpha + \beta) \in \Phi_H$ is equivalent to
$w^{-1}E_{ij}w \in H$, which in turn is just $w(i) \leq h(w(j))$.
\end{proof}

The theorem also simplifies when the operator
is either nilpotent or semisimple.

\section{Open Questions} \label{open questions}

Many questions about Hessenberg varieties remain, some of which
are described here.

\subsection{Geometric properties}  One of the most fundamental unanswered 
questions about the geometry of Hessenberg varieties is:  

\begin{question}
Is every Hessenberg variety pure dimensional?
\end{question}

In every known example, the answer to this is yes.  This also raises the
following.

\begin{question}
What is the dimension of the Hessenberg variety $\H(X,H)$?
\end{question}

The answer is known for various examples, including the Springer fibers
(where it is $\sum_{i=1}^k (i-1) d_i$ if the Jordan blocks have
size $d_1$, $\ldots$, $d_k$ \cite{Sp}) and 
regular nilpotent Hessenberg varieties (namely
$\sum_{i=1}^n (h(i)-i)$ \cite{ST}).  It is unknown in general.

The answer to the next question is known for the
Springer fiber, where it is the dimension of the corresponding 
irreducible representation of the symmetric group (\cite{Sp},
 \cite[3.6.2]{CG}).

\begin{question}
How many components does $\H(X,H)$ have?
\end{question}

This paper has discussed Hessenberg varieties over $GL_n(\C)$.  Hessenberg
varieties are defined for general complex linear algebraic groups (see
\cite{dMPS}), and the same questions can be posed in the general setting.

\begin{question}
How many of these results hold for general $G$?
\end{question}

\subsection{Closure relations} Given a Schubert cell $C_w$, classical
results show that the 
cell $C_x$ lies in its closure if and only if $w$ is a product 
of simple transpositions
$w=s_1 \cdots s_k$ and $x=s_{i_1} \cdots s_{i_{k'}}$ with
$1 \leq i_1 < \cdots < i_{k'} \leq k$ (see \cite[section 2.7]{BL}).  
The closure $\overline{C_w}$
is a Schubert variety, whose geometry and associated combinatorics
has been extensively studied \cite{BL}. 

The cells of a general Hessenberg variety are intersections
with Schubert cells.  However, the closure relations of these
intersections are not in general restrictions of the closure relations of
the full Schubert cells.  

\begin{question}
What are the closure relations for cells in a Hessenberg variety?  
For which $x$ does
$C_x \cap \H(X,H)$ intersect the closure of $C_w \cap \H(X,H)$?
\end{question}

The answer to an apparently simpler question is also unknown.

\begin{question}
If $C_w \cap \H(X,H)$ is nonempty, for which permutations 
$x$ does the flag given by $x$ lie in the closure of $C_w \cap \H(X,H)$?
\end{question}

\subsection{Betti numbers}  The previous results established that the
odd-dimensional Betti numbers for Hessenberg varieties are zero.
They also provide an algorithm to generate tables of the 
even-dimensional Betti numbers, which are available at
\begin{center}{\tt http://www.math.lsa.umich.edu/$\sim$tymoczko}
\end{center}  
The even-dimensional Betti numbers for Hessenberg varieties
$\H(N,H)$  have closed formulae when $N$ is a regular nilpotent operator, 
i.e., $N$
consists of a single Jordan block. 
These Betti numbers are 
both symmetric (namely $b_i = b_{k-i+1}$ for each $i$) and 
unimodal (namely $b_1 \leq b_2 \leq b_3 
\cdots \leq b_{\lceil k/2 \rceil}$) by \cite{ST}.  
Yet most of these varieties are singular.

The even-dimensional Betti numbers for general
Hessenberg varieties need not be symmetric.  Robert MacPherson conjectured
the following, which is true in all known cases.  It is the combinatorial
description of the hard Lefschetz property and has been studied in 
other contexts \cite{Sta}.

\begin{question}
For any Hessenberg variety $\H(X,H)$ the even-dimensional Betti numbers 
are unimodal and satisfy $b_i \leq b_{k-i+1}$ for all $i$ between $1$ and 
$k/2$.
\end{question}

\end{document}